\overfullrule=0pt
\centerline {\bf On the infimum of certain functionals}\par
\bigskip
\bigskip
\centerline {BIAGIO RICCERI}\par
\bigskip
\bigskip
\bigskip
\bigskip
Here and in what follows, $X$ is a real Banach space, $\varphi:X\to {\bf R}$ is a non-zero continuous linear functional and $\psi:X\to {\bf R}$
is a non-constant Lipschitzian functional with Lipschitz constant $L$. \par
\smallskip
From Proposition 2.1 of [2], we kow that, when $L<\|\varphi\|_{X^*}$, the functional $\varphi+\psi$ is unbounded below. When, to the contrary,
$L\geq\|\varphi\|_{X^*}$, this is no longer true. That is, the same functional can be bounded below. The simplest examples are provided  by taking
$\psi(x)=|\varphi(x)|$ or $\psi(x)=\|\varphi\|_{X^*}\|x\|\ .$
\smallskip
The aim of this very short paper is to study the infimum of that functional just when $L=\|\varphi\|_{X^*}\ .$\par
\smallskip
So, from now on, we assume that
$$L=\|\varphi\|_{X^*}\ .$$
\smallskip
Our basic result is as follows:\par
\medskip
THEOREM 1. - {\it Let $[a,b]$ be a closed interval contained in $[-1,1]$ and let $\gamma:[a,b]\to {\bf R}$ be a convex and
lower semicontinuous function.\par
Then, one has
$$\max\left\{\inf_{x\in X}(\varphi(x)+a\psi(x))-\gamma(a),\inf_{x\in X}(\varphi(x)+b\psi(x))-\gamma(b)\right\}=
\inf_{x\in X}\sup_{\lambda\in [a,b]}(\varphi(x)+\lambda\psi(x)-\gamma(\lambda))\ .$$}
\medskip
Our proof of Theorem 1 is based on the use of the following result (Theorem 5.9 and Remark 5.10 of [3]):\par
\medskip
THEOREM A. - {\it Let $S$ be a topological space, $I\subset {\bf R}$ a compact interval and $f:S\times I\to {\bf R}$ a
function which is lower semicontinuous in $X$ and upper semicontinuous and quasi-concave in $I$. Moreover, assume that
there exists a set $D\subseteq I$ dense in $I$ such that, for each $\lambda\in D$ and $r\in {\bf R}$, the set
$$\{x\in S : f(x,\lambda)<r\}$$
is connected.\par
Then, one has
$$\sup_{\lambda\in I}\inf_{x\in S}f(x,\lambda)=\inf_{x\in S}\sup_{\lambda\in I}f(x,\lambda)\ .$$}
\medskip
To be able to use Theorem A, we first have to establish the following result:\par
\medskip
THEOREM 2. - {\it 
For each $\lambda\in ]-1,1[$ and $r\in {\bf R}$,  the set 
$(\varphi+\lambda\psi)^{-1}(]-\infty,r])$ is a retract of $X$.}\par
\smallskip
PROOF. First, consider the multifunction $G:{\bf R}\to 2^X$ defined by
$$G(t)=\varphi^{-1}(]-\infty,t])$$
for all $t\in {\bf R}$. Let us check that
$$d_H(G(t),G(s))\leq{{|t-s|}\over {\|\varphi\|_{X^*}}}\ . \eqno{(1)}$$
for all $t, s\in {\bf R}$, $d_H$ being the usual Hausdorff distance. For instance, assume that $t<s$. Consequently
$$G(t)\subseteq G(s)\ .\eqno{(2)}$$
Now, fix $x\in G(s)\setminus G(t)$. Consequently
$$t<\varphi(x)\leq s\ .$$
In view of the classical formula giving the distance of a point from a closed hyperplane, we have
$$\hbox {\rm dist}(x,G(t))\leq\hbox {\rm dist}(x,\varphi^{-1}(t))
={{\varphi(x)-t}\over {\|\varphi\|_{X^*}}}\leq {{s-t}\over {\|\varphi\|{_X^*}}}\ .$$
So
$$\sup_{x\in G(s)}\hbox {\rm dist}(x,G(t))\leq {{s-t}\over {\|\varphi\|{_X^*}}}$$
which together with $(2)$ gives $(1)$. Now, consider the multifunction $F:X\to 2^X$
defined by
$$F(x)=G(r-\lambda\psi(x))$$
for all $x\in X$. For each $x,y\in X$, we have
$$d_H(F(x),F(y))\leq {{1}\over {\|\varphi\|_{X^*}}}|\lambda||\psi(x)-\psi(y)|\leq |\lambda|\|x-y\|\ .$$
Hence, since $|\lambda|<1$, $F$ is a multivalued contraction with closed and convex values.
Then, in view of [1], the set $Fix(F):=\{x\in X : x\in F(x)\}$ is a retract of $X$. To complete the proof, simply observe
that $Fix(F)=(\varphi+\lambda\psi)^{-1}(]-\infty,r])$.\hfill $\bigtriangleup$\par
\medskip
{\it Proof of Theorem 1}. Consider the function $f:[a,b]\to {\bf R}$ defined by
$$f(x,\lambda)=\varphi(x)+\lambda\psi(x)-\gamma(\lambda)$$
for all $(x,\lambda)\in X\times [a,b]$. Clearly, $f$ is continuous in $X$, while it is upper semicontinuous and concave
in $[a,b]$. Fix $\lambda\in ]a,b[$ and $r\in {\bf R}$. Of course, we have
$$\{x\in X : f(x,\lambda)<r\}=\bigcup_{s<r}\{x\in X : f(x,\lambda)\leq s\}\ .$$
On the other hand, by Theorem 2, the sets of the family $\{\{x\in X : f(x,\lambda)\leq s\}\}_{s<r}$ are connected (being retracts
of $X$) and pairwise non-disjoint. Consequently, the set $\{x\in X : f(x,\lambda)<r\}$ is connected too.
Therefore, we can apply Theorem A. It ensures that
$$\sup_{\lambda\in [a,b]}\inf_{x\in X}(\varphi(x)+\lambda \psi(x)-\gamma(\lambda))=\inf_{x\in X}
\sup_{\lambda\in [a,b]}(\varphi(x)+\lambda \psi(x)-\gamma(\lambda))\ .$$
Now, observe that, since $\inf_{x\in X}(\varphi(x)+\lambda\psi(x))=-\infty$ for all $\lambda\in ]-1,1[$, we have
$$\sup_{\lambda\in [a,b]}\inf_{x\in X}(\varphi(x)+\lambda \psi(x)-\gamma(\lambda))=
\max\left\{\inf_{x\in X}(\varphi(x)+a\psi(x))-\gamma(a),\inf_{x\in X}(\varphi(x)+b\psi(x))-\gamma(b)\right\}$$
and the conclusion follows.\hfill $\bigtriangleup$
\medskip
A consequence of Theorem 1 is as follows:\par
\medskip
THEOREM 3. - {\it Let $[a,b]$ be a closed interval contained in $[-1,1]$ and let $\gamma:[a,b]\to {\bf R}$ be a 
continuous function which is derivable in $]a,b[$. Assume that $\gamma'$ is strictly increasing in $]a,b[$. Set
$$A=\left\{x\in X : \psi(x)\leq \inf_{]a,b[}\gamma'\right\}\ ,$$
$$B=\left\{x\in X : \psi(x)\geq \sup_{]a,b[}\gamma'\right\}$$
and
$$C=\left\{x\in X : \inf_{]a,b[}\gamma'<\psi(x)<\sup_{]a,b[}\gamma'\right\}\ .$$
Finally, denote by $\eta$ the inverse of the function $\gamma'$.\par
Then, one has
$$\max\left\{\inf_{x\in X}(\varphi(x)+a\psi(x))-\gamma(a),\inf_{x\in X}(\varphi(x)+b\psi(x))-\gamma(b)\right\}=$$
$$\min\left \{
\inf_{x\in A}(\varphi(x)+a\psi(x))-\gamma(a), \inf_{x\in B}(\varphi(x)+b\psi(x))-\gamma(b),
\inf_{x\in C}(\varphi(x)+\eta(\psi(x))\psi(x)-\gamma(\eta(\psi(x))))\right \}\ .$$}\par
\medskip
PROOF. Let $f$ be as in the proof of Theorem 1. Fix $x\in X$. Clearly, $f(x,\cdot)$ is concave in $[a,b]$.
Moreover, according to the sign of its
derivative, the function $f(x,\cdot)$ is non-increasing (resp. non-decreasing)
in $[a,b]$ if $x\in A$ (resp. $x\in B$). If $x\in C$, the derivative of $f(x,\cdot)$ vanishes at the point $\eta(\psi(x))$ and so, by concavity,
such a point is the global maximum of $f(x,\cdot)$ in $[a,b]$. Summarizing, we have
$$\sup_{\lambda\in [a,b]}f(x,\lambda)=\cases{\varphi(x)+a\psi(x)-\gamma(a) & if $x\in A$\cr & \cr 
\varphi(x)+b\psi(x)-\gamma(b) & if $x\in B$\cr & \cr \varphi(x)+\eta(\psi(x))\psi(x)-\gamma(\eta(\psi(x))) & if
$x\in C$\ ,\cr}$$ 
and the conclusion clearly follows in view of Theorem 1.\hfill $\bigtriangleup$\par
\medskip
In turn, applying Theorem 3, we obtain the following result:\par
\medskip
THEOREM 4. - {\it We have
$$\max\left\{\inf_{x\in X}(\varphi(x)+\psi(x)),\inf_{x\in X}(\varphi(x)-\psi(x))\right\}=
\inf_{x\in X}(\varphi(x)+|\psi(x)|) \eqno{(3)}$$
and
$$\liminf_{\|x\|\to +\infty}(\varphi(x)+|\psi(x)|)=\inf_{x\in X}(\varphi(x)+|\psi(x)|)\ .\eqno{(4)}$$}\par
\smallskip
PROOF. First, we want to prove that
$$\max\left\{\inf_{x\in X}(\varphi(x)+\psi(x)),\inf_{x\in X}(\varphi(x)-\psi(x))\right\}=
\inf_{x\in X}(\varphi(x)+|\psi(x)|+e^{-|\psi(x)|})\ .\eqno{(5)}$$
Consider the function
$\gamma:[-1,1]\to {\bf R}$ defined by
$$\gamma(\lambda)=\cases{(1-|\lambda|)\log(1-|\lambda|)+|\lambda| & if $|\lambda|<1$\cr & \cr
1 & if $|\lambda|=1$\ .\cr}$$
Clearly, $\gamma$ is continuous in $[-1,1]$, is derivable in $]-1,1[$, $\gamma'$ is strictly increasing and
$\gamma'(]-1,1[)={\bf R}$. Moreover, $\eta$, the inverse of $\gamma'$, is given by
$$\eta(\mu)=\cases{{{|\mu|}\over {\mu}}(1-e^{-|\mu|}) & if $\mu\neq 0$\cr & \cr
0 & if $\mu=0$\ .\cr}$$
So, for each $x\in X\setminus \psi^{-1}(0)$, we have
$$\eta(\psi(x))\psi(x)-\gamma(\eta(\psi(x)))=|\psi(x)|(1-e^{-|\psi(x)|})-(-e^{-|\psi(x)|}|\psi(x)|+1-e^{-|\psi(x)|})=
|\psi(x)|+e^{-|\psi(x)|}-1\ .$$
Clearly, these equalities hold also if $\psi(x)=0$. Consequently, by Theorem 2, after observing that $C=X$, we have
$$\max\left\{\inf_{x\in X}(\varphi(x)-\psi(x))-1,\inf_{x\in X}(\varphi(x)+\psi(x))-1\right\}=
\inf_{x\in X}(\varphi(x)+\eta(\psi(x))\psi(x))-\gamma(\eta(\psi(x))))$$
$$=\inf_{x\in X}(\varphi(x)+|\psi(x)|+e^{-|\psi(x)|})-1$$
which yields $(5)$. Since
$$\max\left\{\inf_{x\in X}(\varphi(x)+\psi(x)),\inf_{x\in X}(\varphi(x)-\psi(x))\right\}\leq
\inf_{x\in X}(\varphi(x)+|\psi(x)|)\leq \inf_{x\in X}(\varphi(x)+|\psi(x)|+e^{-|\psi(x)|})\ ,$$
from $(5)$, we obtain both $(3)$ and
$$\inf_{x\in X}(\varphi(x)+|\psi(x)|)=\inf_{x\in X}(\varphi(x)+|\psi(x)|+e^{-|\psi(x)|})\ .\eqno{(6)}$$
Finally, let us prove $(4)$. Arguing by contradiction, assume that
$$\inf_{x\in X}(\varphi(x)+|\psi(x)|)<\liminf_{\|x\|\to +\infty}(\varphi(x)+|\psi(x)|)\ .$$
Fix $\xi$ satisfying
$$\inf_{x\in X}(\varphi(x)+|\psi(x)|)<\xi<\liminf_{\|x\|\to +\infty}(\varphi(x)+|\psi(x)|)\ .\eqno{(7)}$$
So, there is some $\delta>0$ such that
$$\varphi(x)+|\psi(x)|>\xi \eqno{(8)}$$
for all $x\in X$ satisfying $\|x\|>\delta$. Now, in view of $(6)$, we can fix a sequence $\{x_n\}$ in $X$ such that
$$\lim_{n\to \infty}(\varphi(x_n)+|\psi(x_n)|+e^{-|\psi(x_n)|})=\inf_{x\in X}(\varphi(x)+|\psi(x)|)\ .\eqno{(9)}$$
Clearly
$$\lim_{n\to \infty}(\varphi(x_n)+|\psi(x_n)|)=\inf_{x\in X}(\varphi(x)+|\psi(x)|)\ .\eqno{(10)}$$
In view of $(7)$, there is $\nu\in {\bf N}$ such that
$$\varphi(x_n)+|\psi(x_n)|<\xi$$
for all $n>\nu$. Thus, by $(8)$, we have
$$\sup_{n>\nu}\|x_n\|\leq \delta\ .$$
Then, since $\psi$ is Lipschitzian, the sequence $\{\psi(x_n)\}$ is bounded too. But, $(9)$ and $(10)$ imply that
$$\lim_{n\to \infty}e^{-|\psi(x_n)|}=0$$
which leads to a contradiction. The proof is complete.\hfill $\bigtriangleup$\par
\medskip
We conclude with a consequence of Theorem 3.\par
\medskip
PROPOSITION 1. - {\it Assume that $\psi$ is G\^ateaux differentiable and that both $\varphi$ and
$-\varphi$ do not belong to $\psi'(X)$.\par
Then, for every $r\in {\bf R}$, the functional $x\to \varphi(x)+|\psi(x)-r|$ has no global minima in $X$.}\par
\smallskip
PROOF. Arguing by contradiction, assume that there is $x_0\in X$ such that
$$\varphi(x_0)+|\psi(x_0)-r|=\inf_{x\in X}(\varphi(x)+|\psi(x)-r|)\ .$$
Then, by Theorem 3 (applied to $\psi-r$) , $x_0$ would be a global minimum either of $\varphi+\psi$ or of $\varphi-\psi$. Accordingly,
we would have either $\psi'(x_0)=-\varphi$ or $\psi'(x_0)=\varphi$, contrary our assumption.\hfill $\bigtriangleup$
\medskip
REMARK 1. - Of course, if $\|\psi'(x)\|_{X^*}<L$ for all $x\in X$, then  both $\varphi$ and
$-\varphi$ do not belong to $\psi'(X)$.\par

\vfill\eject

\centerline {\bf References}\par
\bigskip
\bigskip
\noindent
[1]\hskip 5pt B. RICCERI, {\it Une propri\'et\'e topologique de l'ensemble
des points fixes d'une contraction multivoques \`a valeurs convexes},
Rend. Accad. Naz. Lincei, (8) {\bf 81} (1987), 283-286.\par
\smallskip
\noindent
[2]\hskip 5pt B. RICCERI, {\it Further considerations on a variational
property of integral functionals}, J. Optim. Theory Appl., {\bf 106}
(2000), 677-681.\par
\smallskip
\noindent
[3]\hskip 5pt B. RICCERI, 
{\it Nonlinear eigenvalue problems},  
in ``Handbook of Nonconvex Analysis and Applications'' 
D. Y. Gao and D. Motreanu eds., 543-595, International Press, 2010.\par
\bigskip
\bigskip
\bigskip
\bigskip
Department of Mathematics\par
University of Catania\par
Viale A. Doria 6\par
95125 Catania, Italy\par
{\it e-mail address}: ricceri@dmi.unict.it

\bye